# STUDIES IN SIMILARITY

**Christopher J Bradley**

**Abstract:** Three circles define each of the Brocard points of a triangle. If one adds the three circles through a pair of vertices and the orthocentre one has nine circles. It is described how each of the nine centres of these circles lies at the vertex of two triangles producing six triangles. These triangles are each similar to the original triangle, three being directly similar and three indirectly similar. These latter three are mutually in perspective with the circumcentre as perspector.

## 1. Introduction

Given a scalene triangle *ABC* there are precisely six points *J* such that the angles *BJC*, *CJA*, *AJB* have values $\alpha = 180^{\circ} - A$, $\beta = 180^{\circ} - B$, $\gamma = 180^{\circ} - C$ in some order. Three of these points are well known to be the orthocentre *H* and the two Brocard points $\Omega$ and $\Omega'$, which in this article, for reasons that will almost immediately become clear, we denote by *H*+ and *H*-. The other three points are less familiar. We denote them by *aH*, *bH*, *cH* and they lie on the circumference of the orthocentroidal circle (diameter *GH*), where the medians meet this circle (other than at the centroid *G*).

The angles *BJC*, *CJA*, *AJB* at these points are shown in the Table 1 below.

| *J* | $\angle BJC$ | $\angle CJA$ | $\angle AJB$ |
|---|---|---|---|
| *H* | $\alpha$ | $\beta$ | $\gamma$ |
| *H*+ | $\gamma$ | $\alpha$ | $\beta$ |
| *H*- | $\beta$ | $\gamma$ | $\alpha$ |
| *aH* | $\alpha$ | $\gamma$ | $\beta$ |
| *bH* | $\gamma$ | $\beta$ | $\alpha$ |
| *cH* | $\beta$ | $\alpha$ | $\gamma$ |

The notation used for the six points is prompted by Table 1, which reflects the group of symmetries of an equilateral triangle.

Each Brocard point may be constructed by drawing three circles and these circles together with the circles *BHC*, *CHA*, *AHB* produce nine circles. We first show how these nine circles may be used to locate the points *aH*, *bH*, *cH*. We determine the equations of the nine circles and the co-ordinates of *aH*, *bH*, *cH* showing that they have the locations stated above. We also find the co-ordinates of the centres of these nine circles. The nine centres have the following remarkable properties. They form six triangles with each centre appearing as a vertex of two of the triangles.

Each of these triangles is similar to triangle $ABC$, three of them being directly similar and three being indirectly similar (with vertices in the opposite order). Also the three that are indirectly similar are mutually in perspective with $O$, the circumcentre of $ABC$, as perspector.

Given a triangle $ABC$ with orthocentre $H$ then, if we draw a general circle through $H$, its points of intersection $X$, $Y$, $Z$ with the altitudes $AH$, $BH$, $CH$ respectively, determine a triangle $XYZ$ that is indirectly similar to triangle $ABC$. The circles through $H$ are known as Hagge circles [1]. Since the similarity is indirect there is a centre of inverse symmetry $P$ and lines through $P$, which are the axes of inverse symmetry. If $AP$, $BP$, $CP$ meet the circumcircle again at $D$, $E$, $F$ and $XP$, $YP$, $ZP$ meet the Hagge circle again at $U$, $V$, $W$ then triangles $DEF$ and $UVW$ are also similar. There is an obvious connection between the vertices of the two triangles, which is that $U$ is the reflection of $D$ in $BC$, $V$ the reflection of $E$ in $CA$ and $W$ the reflection of $F$ in $AB$.

The reason a Hagge circle carries a triangle that is indirectly similar to $ABC$ by means of the construction mentioned above is that the angles subtended at $H$ by the sides $BC$, $CA$, $AB$ are $\alpha$, $\beta$, $\gamma$ respectively. The proof of this is a straightforward argument using angle chasing.

One might therefore hope, that if we use either of the Brocard points $H+$ or $H-$ and draw a general circle through that Brocard point and find its intersections with the line segments from that Brocard point to the vertices, then these intersections might determine a triangle $ZXY$ or $YZX$, which is indirectly similar to triangle $ABC$. This indeed proves to be the case and again there is a centre and axis of inverse symmetry, see Bradley [2].

If one follows the same procedure with any of the three points $aH$, $bH$, $cH$ we get triangles that turn out to be directly similar to triangle $ABC$. Direct similarities have rather different properties from indirect similarities and these properties are the subject of a comprehensive survey by Wood [3]. We point out how our new circles and triangles fit into Wood's scheme.

## 2. The nine circles and the six triangles

We first consider three circles that determine $H$. If you reflect $O$ in the side $BC$ to obtain the point $aA$ and then draw a circle centre $aA$ of radius $R$ (the circumradius of $ABC$), then this circle is the circle $BHC$. It is also the case that, if $bB$ and $cC$ are the reflections of $O$ in $CA$ and $AB$ respectively, then circles centres $bB$ and $cC$ of radius $R$ are respectively the circles $CHA$ and $AHB$. The common point of these three circles is therefore the orthocentre $H$. In terms of areal co-ordinates, since the co-ordinates of $H$ are given by $(1/(b^2 + c^2 - a^2), 1/(c^2 + a^2 - b^2), 1/(a^2 + b^2 - c^2))$, the equations of these circles are

$BHC$: $$a^2yz + b^2zx + c^2xy - (b^2 + c^2 - a^2)x(x + y + z) = 0, \tag{2.1}$$
$CHA$ $$a^2yz + b^2zx + c^2xy - (c^2 + a^2 - b^2)y(x + y + z) = 0, \tag{2.2}$$
$AHB$ $$a^2yz + b^2zx + c^2xy - (a^2 + b^2 - c^2)z(x + y + z) = 0. \tag{2.3}$$

Using the fact that the centre of a circle is the pole of the line at infinity we find the unnormalised co-ordinates of the centres are
$$aA(- a^2(b^2 + c^2 - a^2), - a^4 - c^4 + a^2b^2 + 2c^2a^2 + b^2c^2, - a^4 - b^4 + 2a^2b^2 + c^2a^2 + b^2c^2), \tag{2.4}$$
$$bB(- b^4 - c^4 + 2b^2c^2 + a^2b^2 + c^2a^2, - b^2(c^2 + a^2 - b^2), - b^4 - a^4 + b^2c^2 + 2a^2b^2 + c^2a^2), \tag{2.5}$$

$$cC(-c^4 - b^4 + c^2a^2 + 2b^2c^2 + a^2b^2, -c^4 - a^4 + 2c^2a^2 + b^2c^2 + a^2b^2, -c^2(a^2 + b^2 - c^2)). \qquad (2.6)$$

Triangle *aAbBcC* is congruent to triangle *ABC* being the image of *ABC* under an 180° rotation about the nine-point centre *N*.

We next consider the three circles defining the Brocard point *H*+, which has co-ordinates $(1/b^2, 1/c^2, 1/a^2)$. These are well known to be (i) the circle through *A* and *B* touching *BC* at *B*, (ii) the circle through *B* and *C* and touching *CA* at *C* and (iii) the circle through *C* and *A* touching *AB* at *A*. We denote their centres by *aB*, *bC*, *cA* respectively. The equations of these circles are:

*AH+B*: $\quad a^2z^2 + (a^2 - b^2)zx - c^2xy = 0, \qquad (2.7)$

*BH+C*: $\quad b^2x^2 + (b^2 - c^2)xy - a^2yz = 0, \qquad (2.8)$

*CH+A*: $\quad c^2y^2 + (c^2 - a^2)yz - b^2zx = 0. \qquad (2.9)$

The unnormalised co-ordinates of their centres are

$aB(2c^2a^2, -a^4 - b^4 + 2a^2b^2 + c^2a^2 + b^2c^2, -c^2(c^2 + a^2 - b^2)), \qquad (2.10)$

$bC(-a^2(a^2 + b^2 - c^2), 2a^2b^2, -b^4 - c^4 + 2b^2c^2 + a^2b^2 + c^2a^2), \qquad (2.11)$

$cA(-c^4 - a^4 + 2c^2a^2 + b^2c^2 + a^2b^2, -b^2(b^2 + c^2 - a^2), 2b^2c^2). \qquad (2.12)$

We next consider the three circles defining the Brocard point *H*-, which has co-ordinates $(1/c^2, 1/a^2, 1/b^2)$. These are well known to be (i) the circle through *A* and *C* touching *BC* at *C*, (ii) the circle through *B* and *A* and touching *CA* at *A* and (iii) the circle through *C* and *B* touching *AB* at *B*. We denote their centres by *aC*, *bA*, *cB* respectively. The equations of these circles are:

*CH-A*: $\quad a^2y^2 - b^2zx + (a^2 - c^2)xy = 0, \qquad (2.13)$

*AH-B*: $\quad b^2z^2 - c^2xy + (b^2 - a^2)yz = 0, \qquad (2.14)$

*BH-C*: $\quad c^2x^2 - a^2yz + (c^2 - b^2)zx = 0. \qquad (2.15)$

The unnormalised co-ordinates of their centres are

$aC(2a^2b^2, -b^2(a^2 + b^2 - c^2), -c^4 - a^4 + a^2b^2 + 2c^2a^2 + b^2c^2), \qquad (2.16)$

$bA(-a^4 - b^4 + b^2c^2 + 2a^2b^2 + c^2a^2, 2b^2c^2, -c^2(b^2 + c^2 - a^2)), \qquad (2.17)$

$cB(-a^2(c^2 + a^2 - b^2), -b^4 - c^4 + c^2a^2 + 2b^2c^2 + a^2b^2, 2c^2a^2). \qquad (2.18)$

The co-ordinates of the nine centres are unnormalised but in each case they sum to the same amount $(a + b + c)(b + c - a)(c + a - b)(a + b - c)$. The nine circles are illustrated in Fig. 1. Also shown in Fig. 1 is the orthocentroidal circle on *GH* as diameter, with equation shown in Bradley and Smith [4] to be

$$(b^2 + c^2 - a^2)x^2 + (c^2 + a^2 - b^2)y^2 + (a^2 + b^2 - c^2)z^2 - a^2yz - b^2zx - c^2xy = 0. \qquad (2.19)$$

The median *AG* meets this circle again at the point *aH* with co-ordinates $aH(a^2, b^2 + c^2 - a^2, b^2 + c^2 - a^2)$. To prove that this is indeed the point *aH* featured in Section 1, one may verify that it lies on each of the circles with centres *aA*, *aB*, *aC* and consequently the sides *BC*, *CA*, *AB* subtend at *aH* the angles $\alpha$, $\gamma$, $\beta$ these being the angles in the segments of those circles. Similarly *bH* and *cH* lie on the orthocentroidal circle and on the medians *BG* and *CG* respectively and have co-ordinates $bH(c^2 + a^2 - b^2, b^2, c^2 + a^2 - b^2)$ and $cH(a^2 + b^2 - c^2, a^2 + b^2 - c^2, c^2)$. The point *bH* lies on the circles with centres *bA*, *bB*, *bC* and *cH* lies on the circles with centres *cA*, *cB*, *cC*.

## 2. Lengths of segments and the similarities

If *ABC* is an acute-angled triangle the distances from the vertices to the orthocentre and to the Brocard points are known, see Shail [5], to be
$AH = 2R \cos A$, $BH = 2R \cos B$, $CH = 2R \cos C$;
$AH+ = 2R_o \sin B/\sin A$, $BH+ = 2R_o \sin C/\sin B$, $CH+ = 2R_o \sin A/ \sin C$;
$AH- = 2R_o \sin C/\sin A$, $BH- = 2R_o \sin A/\sin B$, $CH- = 2R_o \sin B/\sin C$,
where $R_o = ½\, abc/(b^2c^2 + c^2a^2 + a^2b^2)^{1/2}$. The reader may verify that the corresponding formulae for *aH*, *bH*, *cH* are

$A\,aH = 4R_o \cos A$, $B\,aH = 2R_o \sin A/\sin B$, $C\,aH = 2R_o \sin A/ \sin C$,
$A\,bH = 2R_o \sin B/ \sin A$, $B\,bH = 4R_o \cos B$, $C\,bH = 2R_o \sin B/ \sin C$,
$A\,cH = 2R_o \sin C/\sin A$, $B\,cH = 2R_o \sin C/\sin B$, $C\,cH = 4R_o \cos C$,

where $R_o$ is, in each case, the radius of the pedal triangle of the corresponding point. In an obtuse-angled triangle these become signed lengths. For the point *aH* the value of
$R_o = ½bc/(2b^2 + 2c^2 – a^2)^{1/2}$ with a similar formula for the points *bH*, *cH*.

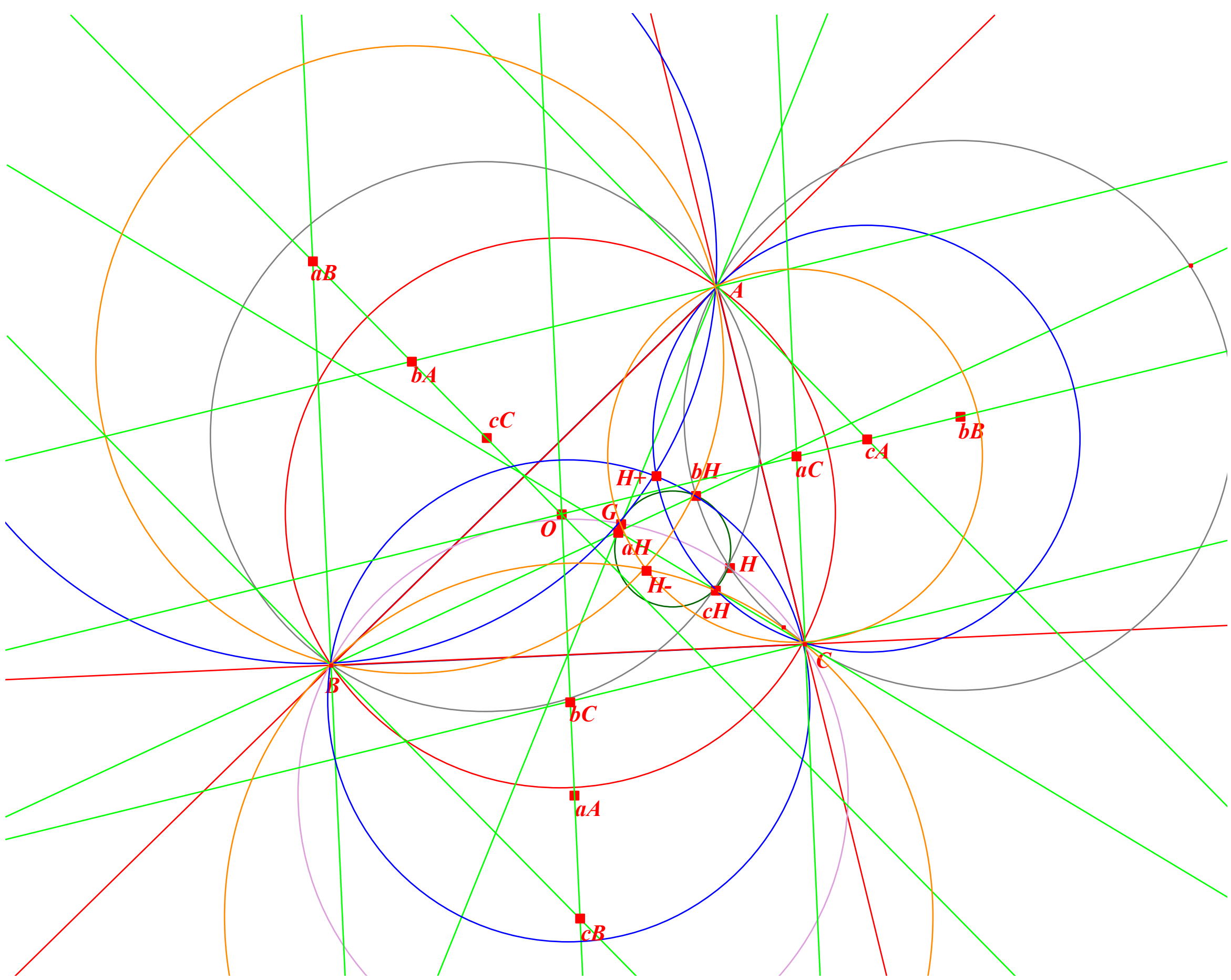

**Fig.1**
**The three circles through each of H, H+ and H- and the orthocentroidal circle**

We next consider the triangle *cA aB bC* and show that it is similar to triangle *ABC*.
The displacement ***bC aB*** $= (a^2(a^2 + b^2 + c^2), - a^4 - b^4 + c^2a^2 + b^2c^2, b^4 - b^2 c^2 - a^2b^2 + 2c^2a^2$
Using the areal metric , see Bradley [6], we find that $(bC\ aB)^2 = ka^2$, where

$$k = (3a^2b^2c^2(a^2 + b^2 + c^2) + 2(b^4c^4 + c^4a^4 + a^4b^4) - a^6(b^2 + c^2) - b^6(c^2 + a^2) - c^6(a^2 + b^2)) \qquad (3.1)$$

all divided by $\{(a + b + c)(b + c - a)(c + a - b)(a + b - c)\}^2$. Similarly $(cA\ bC)^2 = kb^2$ and $(aB\ cA)^2 = kc^2$, from which the result follows. Another short calculation shows that triangle *bA cB aC* is congruent to triangle *cA aB bC*. These two triangles together with *ABC* and *aA bB cC* are all directly similar to one another and fall into two congruent pairs.

We now consider the three triangles *aA aB aC*, *bA bB bC* and *cA cB cC*. Note first that the centres *aA*, *bC*, *cB* all lie on the perpendicular bisector of the side *BC*. Similarly the centres *bB*, *cA*, *aC* all lie on the perpendicular bisector of *CA* and the centres *cC*, *aB*, *bA* all lie on the perpendicular bisector of *AB*. It follows that triangles *aA aB aC*, *bC bA bB*, *cB cC cA* are mutually in perspective with the common perspector *O*.

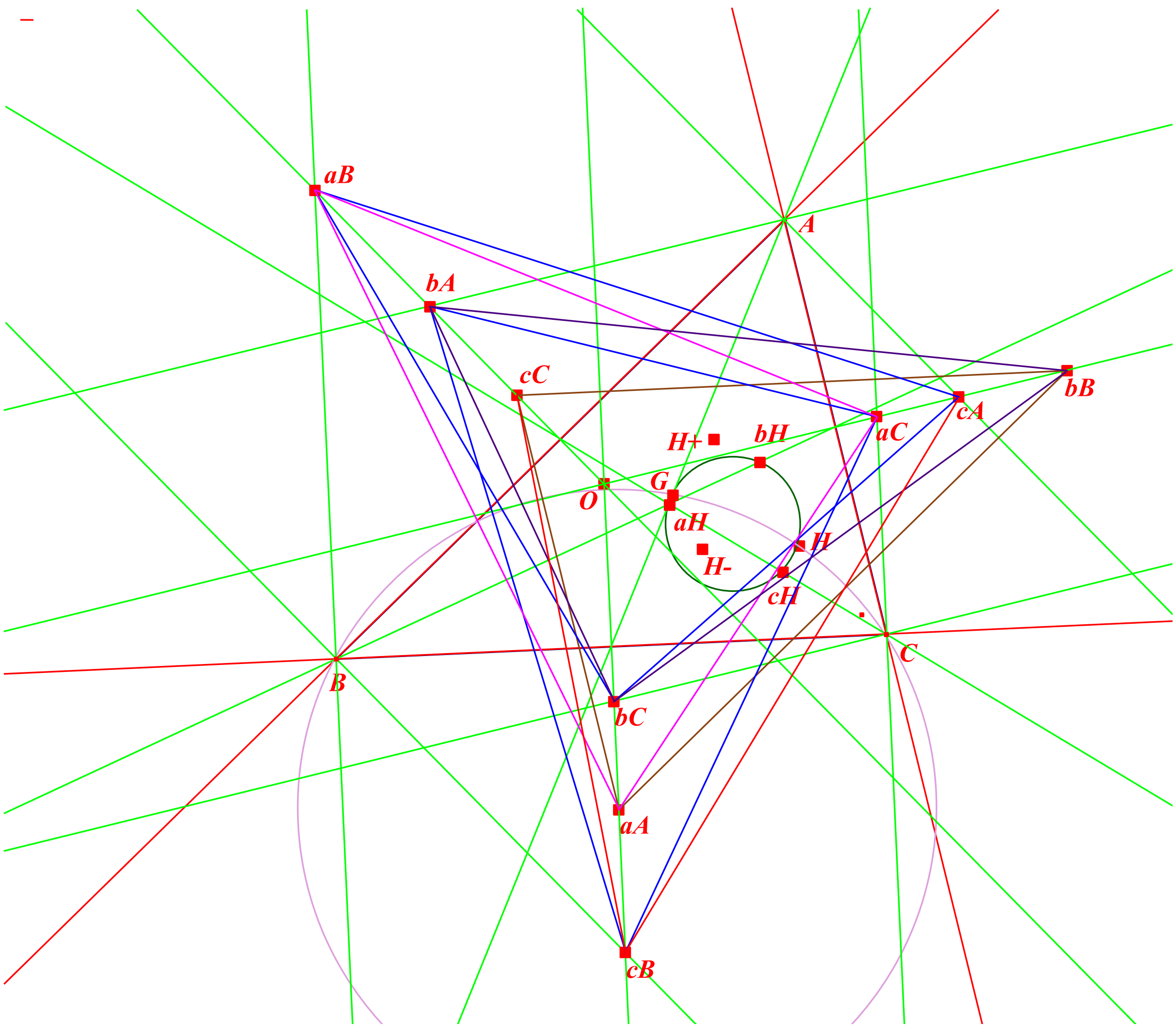

**Fig. 2**
**The nine points each used twice to form the six similar triangles**

The displacement ***aB aC*** = $(2a^2(b^2 - c^2), a^2(a^2 - 3b^2 - c^2), a^2(b^2 + 3c^2 - a^2))$. Using the areal metric we find $(aB\ aC)^2 = la^2$, where

$$l = \{a^2(2b^2 + 2c^2 - a^2)\}/\{(a + b + c)(b + c - a)(c + a - b)(a + b - c)\}. \qquad (3.2)$$

Similar calculations show that $(aC\ aA)^2 = lb^2$ and $(aA\ aB)^2 = lc^2$. It follows that triangle *aA aB aC* is similar to triangle *ABC*.

In similar manner we find $(bB\ bC)^2 = ma^2$, $(bC\ bA)^2 = mb^2$ and $(bA\ bB)^2 = mc^2$, where

$$m = \{b^2(2c^2 + 2a^2 - b^2)\}/\{(a + b + c)(b + c - a)(c + a - b)(a + b - c)\}, \qquad (3.3)$$

showing that triangle *bA bB bC* is also similar to triangle *ABC*.

And again in similar manner we find $(cB\ cC)^2 = na^2$, $(bC\ bA)^2 = nb^2$ and $(bA\ bB)^2 = nc^2$, where

$$n = \{c^2(2a^2 + 2b^2 - c^2)\}/\{(a + b + c)(b + c - a)(c + a - b)(a + b - c)\}, \qquad (3.4)$$

showing that triangle *bA bB bC* is also similar to triangle *ABC*.

All three triangles are directly similar to one another and inversely similar to triangle *ABC*. The areas of the three triangles are in the ratio $l : m : n$. The six triangles are shown in Fig. 2.

Flat 4,
Terrill Court,
12-14 Apsley Road,
BRISTOL BS8 2SP.